\documentclass[letterpaper, 10 pt, conference]{ieeeconf}  

\IEEEoverridecommandlockouts                              

\usepackage{graphics} 
\usepackage{epsfig} 
\usepackage{mathptmx} 
\usepackage{times} 
\usepackage{amsmath} 
\usepackage{amssymb}  

\newtheorem{theorem}{Theorem}

\newtheorem{definition}[theorem]{Definition}
\newtheorem{remark}{Remark}
\newtheorem{assumption}{Assumption}

\title{\LARGE \bf
Partial stabilization of nonlinear systems along a given trajectory}

\author{Victoria Grushkovskaya$^{1,3}$, Iryna Vasylieva$^{1,3}$, and Alexander Zuyev$^{2,3}$
\thanks{$^{1}$Institute of Mathematics, University of Klagenfurt, Universitätsstr. 65--67, 9020 Klagenfurt, Austria
        (Email: {\tt\small viktoriia.grushkovska@aau.at}, {\tt\small iryna.vasylieva@aau.at})}%
\thanks{$^{2}$Max Planck Institute for Dynamics of Complex Technical Systems, Sandtorstra{\ss}e 1, 39106 Magdeburg, Germany (Email: {\tt\small zuyev@mpi-magdeburg.mpg.de})}
\thanks{$^{3}$Institute of Applied Mathematics \& Mechanics, National Academy of Sciences of Ukraine, G.\,Batiuka 19, 84100 Sloviansk, Ukraine}
\thanks{This work is partially supported by the Austrian Science Fund (FWF): DOC 78.}
}

\begin{document}

\maketitle
\thispagestyle{empty}
\pagestyle{empty}

\begin{abstract}
In this paper, the problem of partial stabilization of nonlinear
systems along a given trajectory is considered.
This problem is treated within the framework of stability of a family of sets.
Sufficient conditions for
the asymptotic stability of a one-parameter family of sets using time-dependent
control in the form of trigonometric polynomials are derived. The obtained
results are applied to a model mechanical system.
\end{abstract}

\section{Introduction}

Trajectory tracking is one of the fundamental control problems which has numerous applications in robotics and process engineering. A theoretical justification of tracking properties of control algorithms requires the stability analysis of the tracking error dynamics in a neighborhood of the reference curve.
The stability proof can be straightforwardly achieved, e.g., if the linearized error dynamics is completely controllable. Tracking algorithms, based on the feedback linearization and flatness techniques, are shown to be highly efficient for various engineering models.

For kinematically redundant manipulating robots, the tracking problem can be effectively formulated in terms of a part of the state variables that characterize the control objective.
This analogy creates a connection between pragmatically driven issues in the field of robotics and the notion of partial stability, a concept that was rigorously defined by A.M. Lyapunov and has been extensively explored by many researchers (see, e.g.,~\cite{RO87,Hai95,Zu99,Zu00,Zu01,Che02,Mi02,Sun02,Zu03,Jian04,Mir04,Vo12,Han14,Had15,Zu06,Al17,Qin21} and references therein).
Specifically, the paper~\cite{Hai95} explored the connection between partial stability and full-variable stability in nonholonomic mechanical systems, offering conditions for achieving partial stability. Issues of partial stabilization in the context of Lagrangian systems were investigated in~\cite{Shi00,Ko02}.
Controllers that utilize passivity-based approaches for partial stabilization have been suggested in~\cite{Mir02,Bin11,Wang21}. The issue of achieving partial stabilization in stochastic dynamical systems is addressed, e.g., in the papers~\cite{Raj17,Vor19,ZV19a}.

The issue of achieving partial stabilization within a finite time frame for systems in chained-form and cascade configurations has been investigated, as illustrated in~\cite{Ja14,Chen15,Gol16}. An analysis of more extensive categories of nonlinear control systems can be found in~\cite{Had15,Laff16}. The suggested adequate conditions for partial stability are based on the premise that the system allows for a Lyapunov function, with the time derivative that is negatively definite concerning a specific subset of variables.

In the paper~\cite{rizaldi2018formally}, the motion planning problem for autonomous vehicles is considered within the framework of manoeuvre automata. In order to ensure the safety of paths in complex environments, it is required to estimate the reachable set of each manoeuvre. The proposed motion planning scenario is illustrated by a unicycle mode with two controls corresponding to the angular velocity and the translational acceleration.
The analogy between the equations of motion of nonholonomic systems and underwater vehicles has been pointed out in~\cite{bara}, where driftless control-affine systems have been used to model the kinematics of an autonomous submarine. These equations have been analyzed in the paper~\cite{grushkovskaya2019stabilization} in the context of trajectory tracking problem with oscillating inputs. A survey of recent advances in the motion planning of autonomous underwater vehicles (AUV) is presented in~\cite{panda2020comprehensive}.
A mathematical model of an unmanned surface vehicle (USV) in the form of a nonlinear control-affine system with 6-dimensinal state and 3-dimensional force input is considered in~\cite{sun2021adaptive}. For this dynamical model, a tracking controller is constructed under the assumption that the planar reference trajectory is regular enough and $C^2$-bounded. The stability proof of the tracking algorithm is based on Lyapunov’s direct method.

An important application of partial stability theory originates from planning the motion of robotic systems in task-spaces.
The goal of the latter problem is to steer the output of a nonholonomic system to a neighborhood of the target point. As the number of output variables (which characterize the task space) is usually less than the dimension of the state space, this task fits into the framework of partial stabilization problems. An approach for solving the motion planning problem in task-space is proposed in~\cite{AS21} based on the Campbell--Baker--Hausdorff--Dynkin formula.
The efficiency of this approach has been tested by the unicycle and car models with kinematic control.
Fundamental solutions of the Laplace equation are exploited in~\cite{rousseas2022trajectory} to generate obstacle-free motion of a disk robot in a bounded connected workspace. The control function, corresponding to the robot velocity, is obtained by an appropriate rescaling of the gradient of the potential function. The control scheme is implemented sequentially, and the convergence of the trajectories to the goal is proved. Computational complexity of the proposed control algorithm is estimated by numerical experiments.

While the field of partial stability theory has advanced substantially, contributions to the partial stabilization of underactuated nonlinear control systems remain relatively scarce.
The challenge of partial stabilization persists for general nonholonomic systems due to the difficulty in formulating an appropriate Lyapunov-like function. In~\cite{GZ19a}, practical conditions for partial asymptotic stability were introduced for control-affine systems that exhibit a partially asymptotically stable equilibrium in their averaged form. This paper tackles the issue of devising explicit partially stabilizing feedback mechanisms for nonlinear control-affine systems that comply with a specific Lie algebra rank condition in their vector fields.

This paper presents a novel approach to partial stabilization that significantly advances the state of the art by addressing the challenge of stabilizing along non-feasible curves~-- a task not previously tackled. By conceptualizing this problem through the lens of the stability of sets, we establish a unifies framework that allows for the stabilization of system behaviors in the vicinity of a given trajectory rather than at a fixed point. The introduction of time-varying feedback laws is a crucial point in our construction, ensuring exponential stability across a family of sets proximal to the non-feasible curve.

The rest of this paper is organized as follows. The partial stabilization problem is formulated in Section~II within the framework of a family of sets. The main result (Theorem~1) is presented in Section~III, and its proof is given in the Appendix. Section~IV illustrates our control design scheme for an autonomous underwater vehicle model.

\section{Preliminaries}
\subsection{Notations and definitions}
Consider a nonlinear system of the form
\begin{equation}\label{eq1}
  \dot{x}=f_0(t,x)+\sum_{k=1}^{m}f_k(x)u_k,
\end{equation}
where $x=(x_1,...,x_n)^T \in D\subset {\mathbb R}^n$ is the state vector, $u=(u_1,...,u_m)^T \in \mathbb{R}^m$ is the control,  $m<n$, $f_0:\mathbb R^+\times D\to\mathbb R^n$, and $f_1,\dots,f_m:D\to\mathbb R^n$.
We  represent the state vector as $x=(y^T,z^T)^T$ with $y=(y_1,...,y_{n_1})^T$ and $z=(z_1,...,z_{n_2})^T$, $n_1+n_2= n$,
and assume that $D=D_y\times \mathbb R^{n_2}$, where $D_y\subset {\mathbb R}^{n_1}$ is a domain containing the point $y=0\in\mathbb R^{n_1}$.

We will  consider the problem of stabilization of system~\eqref{eq1} with respect to its $y$-variables.
For this purpose, we introduce some necessary notations and definitions which will be used throughout the paper.

For vector fields $f, g \in C^1 (D; \mathbb{R}^n)$ and a point $x^*\in D,$  we define
the directional derivative
$\mathcal{L}_g f(x^*)=\frac{\partial f(x)}{\partial x}g(x)\bigg|_{x=x^*}$ and the Lie bracket
$[f,g](x^*)=\mathcal{L}_f g(x^*)-\mathcal{L}_g f(x^*)$.
For a time dependent vector field $f \in C^1(\mathbb{R}^+\times D; \mathbb{R}^n)$, the directional derivative at a point $(t^*,x^*)\in \mathbb{R}^+\times D$ is $$\mathcal{L}_g f(t^*,x^*)=\frac{\partial f(t^*,x)}{\partial x}g(x)\bigg|_{{t=t^*},{x=x^*}}.$$

We  say that an  $f:\mathbb{R}^+\times D\to \mathbb{R}^n,$  is:
\begin{itemize}
  \item[-]  \emph{Lipschitz continuous with respect to $x$ uniformly in $t$} in a set $\tilde{D}\subseteq D,$ if there exists an $L>0$ such that $\|f(t,x)-f(t,\tilde{x})\|\le L \|x-\tilde{x}\|$ for all $x,\tilde{x}\in\tilde{D},$ $t\ge0;$
  \item[-] \emph{bounded uniformly in $t$} in a set $\tilde{D}\subseteq D,$ if there exists an $M>0$ such that $\|f(t,x)\|\le M$ for all $x\in\tilde{D},$ $t\ge0.$
\end{itemize}

\begin{definition} \label{def1}
 Given  a time-varying feedback law $u^{\varepsilon}:\mathbb{R}^+\times D\times \mathbb{R}^{n_1} \to \mathbb{R}^m$ depending on a parameter $\varepsilon > 0$ and a vector function $y^{*}: \mathbb{R}^+ \to D_y$, the \textit{$\pi_{\varepsilon}$-solution} of~\eqref{eq1} corresponding to the initial condition $x^0 \in D$ at $t=t_0\ge 0$ and the control $u=u^{\varepsilon}(t,x,y^{*}(t))$ is an absolutely continuous function $x(t)\in D$, defined for $t\in [t_0,+\infty)$, such that $x(t_0)=x^0$ and
\begin{equation}\label{eq2}
\dot{x}=f_0(t,x(t))+\sum_{k=1}^{m}u_k^{\varepsilon}(t,x(t_j),y^*(t_j))f_k(x(t)),
\end{equation}
$$t\in I_j=[t_j,t_{j+1}),\; t_j= t_0 + \varepsilon j\;\text{ for each } j=0,1,2,... \; .
$$
\end{definition}

The concept of $\pi_{\varepsilon}$-solutions has been used, e.g., in \cite{Clar97,ZuSIAM}, and its extension to the case of time-varying control parameters is introduced in \cite{GZ19a}.

\begin{definition}\label{def2}
A one-parametric family of non-empty sets $\{\mathit{Y_t} \}_{t\ge0}$ with $Y_t\subset\mathbb{R}^n$  is called {\em asymptotically stable} for system~\eqref{eq1} with a feedback control of the form $u=u^{\varepsilon}(t,x,y^{*}(t))$, if it is stable and attractive, i.e.:
\begin{itemize}
  \item[$-$] (stability) for every $\Delta>0,$ there exists a $\delta>0$ such that, for every $t_0\ge 0$ and $x^0\in B_{\delta}(\mathit {Y_{t_0}} ),$ the corresponding $\pi_\varepsilon$-solution $x(t)$ with the initial condition $x(t_0)=x^0$ is uniquely defined for $t\ge t_0$  and $x(t)\in\ B_{\Delta}(\mathit {Y_{t}} )$ for all $t\in[t_0;+\infty)$;
  \item[$-$] (attraction) for some $\delta>0$ and for every $\Delta>0$, there exists a $t_1\ge 0$ such that, for any $t_0\ge 0, $ and $x^0\in B_{\delta}(Y_{t_0}),$ the corresponding $\pi_\varepsilon$-solution $x(t)$ with the initial condition $x(t_0)=x^0$ satisfies the property
   $ x(t)\in B_{\Delta}(\mathit {Y_{t}} )\text{ for all } t\in[t_0+t_1,\infty).$
\end{itemize}
\end{definition}

The stability concept for families of sets is described, e.g., in~\cite{La02,GZ18,GZ19a}.

\subsection{Problem statement}

Let $D_y \subset \mathbb{R}^{n_1}$ be a non-empty domain, $y^*(t)$ be a curve in $D_y,$ $y^*\in C(\mathbb{R}^+; D_y). $ In this paper, we consider the following problem of stabilizing the $y$-variables of system~\eqref{eq1} along the curve $y^*(t):$

{\bf Problem~1.}
\emph{Given a curve $y^*\in C(\mathbb{R}^+; D_y)$ and a number $p>0,$ the goal is to find a control $u^{\varepsilon}(t,x,y^*)$ such that the family of sets $\{\mathit{Y_t^p}\}_{t\ge0}$ with
\begin{equation}\label{sets}
\begin{aligned}
Y_t^p=\{x=(y^T,z^T)^T\in\mathbb{R}^{n_1}:\|y^*(t)-y\|<p, z \in \mathbb{R}^{n_2}\}
\end{aligned}
\end{equation}
 is asymptotically stable for the closed-loop system~\eqref{eq1} with $u = u^{\varepsilon}(t,x,y^*)$ in the sense of Definitions 1 and 2.}

In the sequel, by a neighborhood of a set $Y_t^p,$ $t\ge 0,$ we mean the set $B_{\delta}(Y_t^p)=\{x\in\mathbb{R}^n: \|y-y^*\|\le p+ \delta, z\in \mathbb{R}^{n_2}\}$.
We assume that $p$ is small enough to guarantee that $B_p (y^*(t))\subset D_y$ for all $t\ge 0$.

Note that the partial stabilization problem for the case of static $y^*(t) \equiv const$ is considered in \cite{GZ20_ECC} under appropriate controlability rank condition. In the paper \cite{GZ19a}, the problem of stabilizing the trajectories of a nonholonomic system to a reference curve in $\mathbb{R}^n$ is considered. Up to our best knowledge, the problem of \emph{partial stabilization to a curve} is considered here for the first time for underactuated nonlinear systems.

For clarity of presentation, we rewrite system~\eqref{eq1} as
\begin{equation}\label{eq3}
    \dot{y}=g_0(t,x)+\sum_{k=1}^{m}g_k(x)u_k,\;\;    \dot{z}=h_0(t,x)+\sum_{k=1}^{m}h_k(x)u_k,
\end{equation}
where $g_k:\mathbb{R}^n\to\mathbb{R}^{n_1}$, $h_k:\mathbb{R}^n\to\mathbb{R}^{n_2}$, $g_0:\mathbb R^+\times\mathbb{R}^n\to\mathbb{R}^{n_1}$ and $h_0:\mathbb R^+\times\mathbb{R}^n\to\mathbb{R}^{n_2}$  are such that the vector fields of  system~\eqref{eq1} are represented as
$$f_0(t,x) = \left(\begin{array}{c}
                      g_0(t,x)\\
                      h_0(t,x)
                    \end{array}\right), f_k(x)=\left(\begin{array}{c}
                                                            g_k(x) \\
                                                            h_k(x)
                                                          \end{array}\right),  k=\overline{1,m}. $$
\section{Main result}

In this section, we consider the class of systems~\eqref{eq3}, whose control vector fields $g_k\in C^1(\mathbb R^n;\mathbb R^{n_1})$ satisfy the following rank condition
for all $x\in D$:
\begin{equation}\label{eq4}
  \text{ span}\left\{\big(g_i(x)\big)_{i\in S_1},\big(I_{n_1\times n}[f_{i_1},f_{i_2}](x)\big)_{(i_1,i_2)\in S_2}\right\}=\mathbb{R}^{n_1},\end{equation}
  where $S_1$, $S_2$ are some sets of indices $S_1\subseteq\{1,2,...,m\},$ $S_2\subseteq\{1,2,...,m\}^2$ such that $|S_1|+|S_2|=n_1$.

This assumption represents a relaxation of the controllability rank condition that the vector fields of system~\eqref{eq1} with their Lie brackets span the whole tangent space:
$$
  \text{ span}\left\{\big(f_i(x)\big)_{i\in \widetilde S_1},\big([f_{i_1},f_{i_2}](x)\big)_{(i_1,i_2)\in \widetilde S_2}\right\}=\mathbb{R}^{n}
  $$
 at each $x\in D$ with some $\widetilde S_1\subseteq\{1,2,...,m\},$ $\widetilde S_2\subseteq\{1,2,...,m\}^2$, $|\widetilde S_1|+|\widetilde S_2|=n$.
For the partial stabilization problems, the latter requirement  can be replaced with relaxed condition~\eqref{eq4}. Thus, we take into account only the first $n_1$ coordinates of $f_i(x)$ and $[f_{i_1},f_{i_2}](x)$, i.e. we exploit the vector fields $g_i(x)$ and $I_{n_1\times n}[f_{i_1},f_{i_2}](x)=\mathcal L_{f_{i_1}}(x)g_{i_2}(x)-\mathcal L_{f_{i_2}}(x)g_{i_1}(x)$. Consequently, a smaller set of vector fields is needed to satisfy the stabilizability condition, which simplifies the control design. Condition~\eqref{eq4}  has been proposed in \cite{GZ20_ECC} for the case $y^*(t) \equiv const,$ and we exploit it here for solving Problem 1.

In order to stabilize the $y$-variables of system~\eqref{eq1} along a curve $y^*(t),$ we will use a time-varying feedback control of the form
\begin{equation}
\label{eq5}
\begin{aligned}
 u_k^\varepsilon(t,x,y^*)=&\sum\limits_{i\in S_1}\phi_i^k(t,x,y^*)\\
 &+\frac{1}{\sqrt{\varepsilon}}\sum\limits_{(i_1i_2)\in S_2}\phi_{i_1,i_2}^k(t,x,y^*),k=\overline{1,m},
\end{aligned}
\end{equation}
  where
  \begin{equation*}
  \begin{aligned}
    &\phi_i^k(t,x,y^*)=\delta_{ki} a_i(x,y^*),\\
 & \phi_{i_1,i_2}^k(t,x,y^*)=2\sqrt{\pi\kappa_{i_1i_2}|a_{i_1i_2}(x,y^*)|}\left( \delta_{ki_1}\cos\left(\frac{2\pi\kappa_{i_1i_2}t}{\varepsilon}\right)\right.\\
 & \qquad\qquad\left.+\delta_{ki_2}{\rm{sign}}(a_{i_1
  i_2}(x,y^*))\sin\left(\frac{2\pi\kappa_{i_1i_2}t}{\varepsilon}\right)\right).
\end{aligned}
\end{equation*}

Here, $\varepsilon>0$ is a small parameter, $\kappa_{i_1i_2}\in\mathbb{N}$ are pairwise distinct numbers, $\delta_{ij}$ is the Kronecker delta, and $\left( (a_i(x,y^*))_{i\in S_1},(a_{i_1i_2}(x,y^*))_{(i_1i_2)\in S_2}\right)^T=a(x,y^*),$ where
\begin{equation}\label{eq6}
  a(x,y^*)=-\alpha \mathcal{F}^{-1}(x)(y-y^*),\quad \alpha>0,
\end{equation}
with $\mathcal{F}^{-1}(x)$ denoting the inverse for $n_1\times n_1$ matrix
 $$\mathcal{F}(x)=\left( (g_i(x))_{i\in S_1},(I_{n_1\times n}[f_{i_1},f_{i_2}](x))_{(i_1i_2)\in S_2}\right).$$

 Obviously, the matrix $\mathcal{F}(x)$ is nonsingular in $D$ because of condition~\eqref{eq4}.

 Let us mention that  controllers of the form~\eqref{eq5}-\eqref{eq6} has been used, e.g. in \cite{ZuSIAM,GZ19a,GZ20_ECC}. In this paper, we adopt the control design from the above mentioned papers to solve Problem 1.
  Before formulating the main result of this section, we introduce several assumptions on the vector field of system~\eqref{eq3} and the curve $y^*(t)$.
 \begin{assumption}
   We suppose that the following properties hold in $D=D_y\times \mathbb{R}^{n_2}.$
\begin{itemize}
                        \item[A1.1)] The functions $f_k\in C^1(D;\mathbb R^n)$, $k=\overline{1,m}$, satisfy the rank condition~\eqref{eq4}. Moreover,    $g_k \in C^2(D;\mathbb{R}^{n_1})$,  $g_0 \in C^1(\mathbb R^+\times D; \mathbb{R}^{n_1}),$ and $h_0 \in C(\mathbb R^+\times D; \mathbb{R}^{n_1})$.
                     \item[A1.2)] For any compact set $\widetilde{D_y}\subset D_y,$ for all $k_1, k_2,k_3\in\overline{1,m},$ $j_1,j_2\in\overline{0,m}$,
                     \begin{itemize}
                       \item the functions $f_{k_1},$ $\mathcal L_{f_{k_2}}g_{k_1},$ $\mathcal L_{f_{k_1}}\mathcal L_{f_{k_2}}g_{k_1}$ are bounded in $\widetilde{D}=\widetilde{D_y}\times \mathbb{R}^{n_2};$
                       \item the functions $f_{k_1}$ are Lipschits continuous in $\widetilde{D}$;
                     \item the functions $f_{0},$ $\mathcal L_{f_{j_1}}g_{j_2},$ $\mathcal L_{f_{0}}L_{f_{k_2}}g_{k_1}$ and $\frac{\partial g_0}{\partial t}$ are bounded uniformly in $t$ in $\widetilde{D};$
                    \item the function $f_{0}$ is Lipschits continuous with respect to $x$ uniformly in $t$ in $\widetilde{D}.$
                     \end{itemize}
                    \item[A1.3)] The function $y^*:\mathbb{R}^+\to D_y$ is Lipschitz continuous.
                   \end{itemize}
 \end{assumption}

 The following result shows that the family of controls~\eqref{eq5}-\eqref{eq6} solves  Problem 1 for system~\eqref{eq3} under Assumption~1.

{\bf Theorem~1.}
  \textit{Let Assumption 1 be satisfied for system~\eqref{eq3} and a curve  $y^* \in C(\mathbb{R}^+; D_y)$, and let $p,\delta >0$ be arbitrary numbers such that $B_{\delta}(\mathit{Y_t^p}) \subset D$ for all $t\ge 0,$ where the sets $\mathit{Y_t^p}$ are defined in~\eqref{sets}.}

  \textit{Then there exists an $\overline{\varepsilon} >0 $ such that, for any $\varepsilon \in (0,\overline{\varepsilon}],$ the family of sets $\left\{\mathit{Y_t^p}\right\}_{t\ge 0}$  is asymptotically (and even exponentially) stable for system~\eqref{eq3} with the controls $u_k = u_k^{\varepsilon}(t,x,y^*)$ defined by~\eqref{eq5} and the initial conditions $x(0) = x^0 \in B_{\delta}(\mathit{Y_t^p}).$}

The proof of this theorem is presented in the Appendix.

\begin{remark}
Unlike the paper \cite{GZ20_ECC}, we do not require the $z$-extendability of solutions to system~\eqref{eq3}, which is instead guaranteed by Assumptions A1.1)--A1.2).
 However, if it holds that $z(t)$-variables of the solutions of system~\eqref{eq3} belong to some set $D_2\subset \mathbb R^{n_2}$ whenever the corresponding part $y(t)$ is in $D_1$, than we can take $\widetilde D= D_y\times D_2$ in A1.2). If, additionally, the functions $g_0,g_k$ are bounded uniformly in $t$  in $\widetilde D$, then the boundedness and Lipschitz continuity properties of the functions $h_0,h_k$ are not required.
This can be easily seen from the proof of Theorem~1.
\end{remark}

\begin{remark}
With the use of control formulas from~\cite{GZ18}, the obtained result can be easily extended to systems whose vector fields satisfy the controllability rank condition with first- and second-order Lie brackets.
\end{remark}

\section{Case study: an autonomous underwater vehicle model}
Consider the equations of motion of an autonomous underwater vehicle with four independent controls:
\begin{equation}\label{un_1}
  \begin{array}{lcl}
   \dot x_1=\cos(x_5)\cos(x_6)v,\\
   \dot x_2=\cos(x_5)\sin(x_6)v,\\
   \dot x_3=-\sin(x_5)v,\\
   \dot x_4= \omega_1 + \omega_2\sin(x_4)\tan(x_5)+\omega_3\cos(x_4)\tan(x_5),\\
   \dot x_5=\omega_2\cos(x_4)-\omega_3\sin(x_4),\\
   \dot x_6=\omega_2\sin(x_4)\sec(x_5)+\omega_3\cos(x_4)\sec(x_5).
  \end{array}
\end{equation}
Here, $(x_1,x_2,x_3)$ denote the position of the center of mass, $(x_4,x_5,x_6)$ describe the vehicle orientation (Euler angles), $v$ is the  translational velocity along the $Ox_1$ axis, and $\omega_1, \omega_2, \omega_3$ are the angular velocity components.
Such equations of motion have been presented, e.g., in~\cite{bara}. The stabilization problem for system~\eqref{un_1} by means of oscillating control is considered in~\cite{GZ19a}.
In this section, we consider the problem of stabilizing the  $(x_1,x_2,x_3)$ coordinates of system~\eqref{un_1} by three controls $v$, $\omega_2$, and $\omega_3$, so we assume that the first
component of the angular velocity cannot be controlled.
 Let us denote $v=u_1, \omega_2 = u_2, \omega_3 = u_3,$ $y=(x_1,x_2,x_3)^T$, $z=(x_4,x_5,x_6)^T$, $x=(y^T,z^T)^T$, and rewrite system~\eqref{un_1} in the form~\eqref{eq3}:
\begin{equation*}
\begin{aligned}
  & \dot{y} = \sum_{k=1}^{3} u_k g_k(x),\;\dot z= h_0(t,x)+\sum_{k=1}^{3} u_k h_k(x),
\end{aligned}
 \end{equation*}
 where
\begin{equation*}
\begin{aligned}
h_0(t,x)& = (\omega_1(t),0,0)^T,\\
g_1(x)& = (\cos x_5 \cos x_6, \cos x_5 \sin x_6, -\sin x_5)^T,\\
g_2(x)& = g_3(x) =h_1(x) = (0, 0, 0)^T,\\
h_2(x)& = (\sin x_4 \tan x_5, \cos x_4, \sin x_4\sec x_5)^T,\\
h_3(x)& = (\cos x_4 \tan x_5, -\sin x_4, \cos x_4\sec x_5)^T.
\end{aligned}
\end{equation*}
The rank condition~\eqref{eq4} is satisfied in $D=\{x\in\mathbb R^6:-\frac{\pi}{2}<x_5<\frac{\pi}{2}\}$ with
 $S_1 = \{1\},$ $S_2 = \{(1,2),(1,3)\}.$ Indeed, it is easy to check that the matrix $\mathcal{F}(x)$ is nonsingular in $D$ with $\det \mathcal{F}(x)\equiv 1$:
 $$
 \mathcal F(x)=\left(g_1(x)\quad I_{3\times 6}[f_1,f_2](x)\quad I_{3\times 6}[f_1,f_3](x)\right),
 $$
 where $$
 \begin{aligned}
 I_{3\times 6} &[f_1,f_2](x)=
\big( \cos x_4\sin x_5\cos x_6+\sin x_4\sin x_6,\\
&\cos x_4\sin x_5\sin x_6-\sin x_4\cos x_6,\, \cos x_4 \cos x_5
 \big)^T,\\
  I_{3\times 6} &[f_1,f_3](x)=
\big( -\sin x_4\sin x_5\cos x_6+\cos x_4\sin x_6,\\
&-\sin x_4\sin x_5\sin x_6-\cos x_4\cos x_6,\, -\sin x_4\cos x_5
 \big)^T.
  \end{aligned}
 $$
According to formulas~\eqref{eq5}, we define the controls as $u_k =  u_k^{\varepsilon}(t,x,y^*)$, so that
\begin{equation}\label{un_3}
\begin{aligned}
u_1=&a_1(x,y^*)+ \sqrt{\frac{4 \pi \kappa_{12} |a_{12}(x,y^*)|}{\varepsilon}}\cos\frac{2 \pi \kappa_{12} t}{\varepsilon}\\
&+\sqrt{\frac{4 \pi \kappa_{13} |a_{13}(x,y^*)|}{\varepsilon}}\cos\frac{2 \pi \kappa_{13} t}{\varepsilon},\\
 u_2=& {\rm{sign}}(a_{12}(x,y^*))\sqrt{\frac{4 \pi \kappa_{12} |a_{12}(x,y^*)|}{\varepsilon}}\sin\frac{2 \pi \kappa_{12} t}{\varepsilon},\\
 u_3=&{\rm{sign}}(a_{13}(x,y^*))\sqrt{\frac{4 \pi \kappa_{13} |a_{13}(x,y^*)|}{\varepsilon}}\sin\frac{2 \pi \kappa_{13} t}{\varepsilon},
\end{aligned}
\end{equation}
where
$a(x,y^*) = -\alpha \mathcal{F}^{-1}(x)(y-y^*).$

For numerical simulations, we choose ${y^*(t) = (0.2 t\cos(0.2 t),0.2 t\sin(0.2 t),0.2 t)^T}$, $\omega_1(t) = 0.25\cos(t)$ and put $\varepsilon = 0.1,$ $\alpha = 15$. Fig.~1 illustrates the behavior of system~\eqref{un_1} with control~\eqref{un_3} and the initial condition ${x^0=(0,0,0,\frac{\pi}{4},\frac{\pi}{4},\frac{\pi}{4})^T}$.
\begin{figure}[h!]
\center{
\includegraphics[scale=0.5]{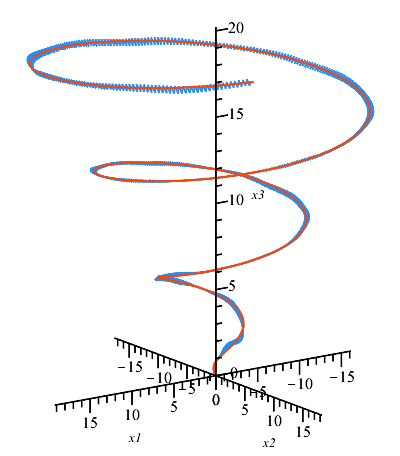}
}
\caption{The blue graph illustrates the behavior of $(x_1,x_2,x_3)$-components of the solution of system~\eqref{un_1}--\eqref{un_3},
and the red curve is $y^*(t)$.}
\end{figure}


\section{Conclusion}
The presented case study demonstrates that our approach is applicable
to the class of underactuated control-affine systems adhering to a specific Lie algebra rank condition, thereby encompassing essentially nonlinear dynamical behavior. On one hand, our method extends the paradigm of partial stabilization to encompass curve-following behaviors; on the other hand, it generalizes our earlier results by encompassing systems with non-zero drift and situations where the reference curve is defined in a lower dimensional subspace.
The outcome of this work is oriented towards robotics, where there is a compelling need to stabilize only a portion of a system's states, enhancing the control and maneuverability of robotic platforms across diverse and potentially unpredictable environments.

\appendix

The proof of Theorem~1 combines and extends the techniques introduced in the papers~\cite{GZ19a,GZ20_ECC,GZ23}. Note that the results of those papers cannot be directly applied because of more general assumptions. In particular, we do not require the $z$-extendability of solutions.

{\em Proof of Theorem~1.}  Given a $p>0,$ let us fix $\delta,$ $\delta^{'}$ such that $0<\delta<\delta^{'}$ and $B_{\delta^{'}}(y^*(t))\subset D_y$ for all $t\ge 0.$ Denote $D^{'} = B_{\delta^{'}}(y^*(t))\times \mathbb{R}^{n_2}.$
   From  assumption A1.2), there exist positive constants $M_g,$ $M_h,$ $M_{g_0},$ $M_{h_0},$ $L_g,$ $L_h,$ $L_{g_0},$ $L_{h_0},$ $M_{g_2},$ $M_{g_{20}},$ $L_{g_{20}},$ $M_{g_3},$ $M_{g_{30}}$  such that for all $x, \tilde{x} \in D^{'},$ $t, \tilde{t} \ge 0,$ and $k_1, k_2, k_3 \in\{1,2,3,...,m\},$
   \begin{equation*}
                     \begin{aligned}
                     &   \|y^*(t){-}y^*(\tilde{t})\| \le L^*|t{-}\tilde{t}|,\\
                     &  \|g_0(t,x)\| \le M_{g_0}, \, \|g_{k}(x)\|\le M_g,\, \left\|\frac{\partial g_0(t,x)}{\partial t} \right\|\le L_{g_{20}},\\
                     & \|h_0(t,x)\| \le M_{h_0}, \, \|h_{k}(x)\|\le M_h,\,  \|\mathcal L_{f_{k_2}}g_{k_1}(x)\| \le M_{g_2},\\
                     & \|g_0(t,x){-}g_0(t,\tilde{x})\|\le L_{g_0} \|x{-}\tilde{x}\|,  \\
                     & \|h_0(t,x){-}h_0(t,\tilde{x})\|\le L_{h_0} \|x{-}\tilde{x}\|, \\
                     &\|g_k(x){-}g_k(\tilde{x})\|\le L_{g} \|x{-}\tilde{x}\|,\|h_k(x){-}h_k(\tilde{x})\|\le L_{h} \|x{-}\tilde{x}\|,\\
                     &    \max\{\|\mathcal L_{f_0}g_0(t,x)\|, \|\mathcal L_{f_{k}}g_0(t,x)\|,\,\|\mathcal L_{f_0}g_{k}(t,x)\|\} \le M_{g_{20}},\\
                     &     \|\mathcal L_{f_{k_3}}\mathcal L_{f_{k_2}}g_{k_1}(x)\|\le M_{g_3},\,    \|\mathcal L_{f_{0}}\mathcal L_{f_{k_2}}g_{k_1}(t,x)\|\le M_{g_{30}}.
                     \end{aligned}
   \end{equation*}

Furthermore,  assumption A1.1) implies the existence of a $\mu >0$ such that $\|\mathcal{F}^{-1}(x)\|\le \mu$ for all $x\in D',$ where $\mathcal{F}^{-1}(x)$ is the inverse matrix for $\mathcal{F}(x).$

 Let $x^0=(y^{0^T},z^{0^T})^T\in B_{\delta}(\mathit{Y_{0}^p}),$ and denote
 $${U^{\varepsilon}=\max\limits_{0\le t\le \varepsilon} \sum_{k=1}^{m}|u_k^\varepsilon(t,x^0,y^*_0)|.}$$
  For the simplicity and without loss of generality, we put $t_0=0$.
 Using~\eqref{eq6} and H\"{o}lder's inequality, one can show that, for any $x^0\in B_{\delta}(\mathit{Y_{0}^p}) ,  $
\begin{equation}\label{eq15}
  U^{\varepsilon}\le c_{1}\|y^0-y_0^*\|+\frac{c_{2}}{\sqrt{\varepsilon}}\sqrt{\|y^0-y_0^*\|}\le c_u \sqrt{\frac{\|y^0-y_0^*\|}{\varepsilon}},
\end{equation}
where $c_1{=}\sqrt{|S_1|}\alpha\mu,$  ${c_2{=}2\sqrt{2\alpha\pi\mu}\left( \sum_{(j_1j_2)\in S_2}(\kappa_{j_1j_2})^{\frac{2}{3}}\right)^\frac{3}{4}}$, and  $c_u=c_1\sqrt{\varepsilon(p+\delta)}+c_2$.

The first step of the proof is to show that all solutions of system~\eqref{eq3} with initial conditions in $B_{\delta}(\mathit{Y_{0}^p}) $  are well defined in $D'$ on the time interval $[0,\varepsilon]$ with some small enough ${\varepsilon>0.}$ Using the integral representation of the $y$-component of the solutions of system~\eqref{eq3} with $x^0\in B_{\delta}(\mathit{Y_{0}^p})$, we get
\begin{equation*}
 \begin{aligned}
\|y(t){-}y^0\|=&\Big{\|}\int\limits_{0}^{t}g_0(s,x(s))ds{+}\sum_{k=1}^{m}\int\limits_{0}^{t} g_k(x(s))u_k(s)ds\Big{\|}\\
\le &\int\limits_{0}^{t}\|g_0(s,x^0)\|ds{+}\sum_{k=1}^{m}\|g_k(x^0)\|\int\limits_{0}^{t}|u_k(s)|ds\\
 &{+}\int\limits_{0}^{t}\|g_0(s,x(s)){-}g_0(s,x^0)\|\\
 &{+}\sum_{k=1}^{m}\|g_k(x(s)){-}g_k(x^0)\||u_k(s)|ds\\
\le & (L_0{+}L_g U^{\varepsilon}) \int\limits_{0}^{t} \Big(\|y(s){-}y^0\| {+} \|z(s){-}z^0\| \Big)ds\\
&{+}(M_{g_0}{+}M_g U^{\varepsilon})t.
 \end{aligned}
 \end{equation*}
Similarly,
\begin{equation*}
\begin{aligned}
    \|z(t){-}z^0\|&\le (L_{h_0}{+}U^{\varepsilon}L_h)\int\limits_{0}^{t}\left(\|y(s){-}y^0\|{+}\|z(s){-}z^0\|\right) ds\\
   &{+}(M_{h_0}{+}U^{\varepsilon}M_h)t.
\end{aligned}
\end{equation*}
Applying Gr\"{o}nwall--Bellman inequality to the both estimates, we obtain:
\begin{equation*}
\begin{aligned}
   \|y(t)-y^0\|\le& e^{(L_0{+}L_gU^{\varepsilon})t}\Big((M_{g_0}{+}U^{\varepsilon}M_g)t\\
   & {+}(L_{g_0}{+}U^{\varepsilon}L_g)\int\limits_{0}^{t}\|z(s)-z^0\|ds\Big),
\end{aligned}
\end{equation*}
\begin{equation*}
\begin{aligned}
    \|z(t)-z^0\|\le&e^{(L_{h_0}{+}U^{\varepsilon}L_h)t}\Big( (M_{h_0}{+}U^{\varepsilon}M_h)t\\
    &{+}(L_{h_0}{+}U^{\varepsilon}L_h)\int\limits_{0}^{t}\|y(s)-y^0\|ds\Big).
\end{aligned}
\end{equation*}

Thus, for any $t\in[0,\varepsilon],$ $x^0\in B_{\delta}(Y_0^p),$
\begin{equation}\label{eq7}\begin{aligned}
  \|y(t)-y^0\|\le e^{c_g\sqrt\varepsilon}\Big(M_{g_0}\varepsilon &+ M_g c_u \sqrt{\varepsilon \|y^0-y^*_0\|}\\
  &+\frac{c_g}{\sqrt\varepsilon}\int\limits_{0}^{t}\|z(s)-z^0\|ds\Big),
\end{aligned}\end{equation}
\begin{equation}\label{eq8}\begin{aligned}
  \|z(t)-z^0\|\le e^{c_h\sqrt\varepsilon}\Big(M_{h_0}\varepsilon &+ M_h c_u \sqrt{\varepsilon \|y^0-y^*_0\|}\\
  &+\frac{c_h}{\sqrt\varepsilon}\int\limits_{0}^{t}\|y(s)-y^0\|ds\Big),
\end{aligned}\end{equation}
where
$$c_g {=} L_{g_0}\sqrt\varepsilon+L_g c_u \sqrt{p+\delta}, \,c_h{ =} L_{h_0}\sqrt\varepsilon+L_h c_u \sqrt{p+\delta}.$$
Substituting~\eqref{eq8} into~\eqref{eq7}, we get:
\begin{equation*}
\begin{aligned}\|y(t)-y^0\|&\le
e^{c_g\sqrt\varepsilon}\Big( M_{g_0}\varepsilon+M_g c_u \sqrt{\varepsilon\|y^0-y_0^*\|}
\\&+c_ge^{c_h\sqrt\varepsilon}\Big(M_{h_0}\varepsilon^{3/2}+M_hc_u\sqrt{\varepsilon\|y^0-y_0^*\|}\\
&+\frac{c_h}{\varepsilon}\int\limits_{0}^{t}\int\limits_{0}^{s}\|y(p)-y^0\|dpds\Big)\Big).
\end{aligned}
\end{equation*}
Then integration by part in the last term of the above estimate yields:
\begin{equation*}\begin{aligned}
\|y(t)-y^0\|\le& e^{c_g\sqrt\varepsilon}\Big(\varepsilon(M_{g_0}+c_g e^{c_h\sqrt\varepsilon}M_{h_0}\sqrt\varepsilon)\\
&+c_u \sqrt{\varepsilon\|y^0-y_0^*\|}(M_g+c_ge^{c_h\sqrt{\varepsilon}}M_h\sqrt\varepsilon)\Big)\\
&+c_gc_h e^{(c_g+c_h)\sqrt\varepsilon}\int\limits_{0}^{t}\|y(s)-y^0\|ds .
\end{aligned}\end{equation*}

Applying again Gr\"{o}nwall--Bellman inequality, we conclude that, for any $\varepsilon>0$ and for all $t\in [0,\varepsilon],$
\begin{equation}\label{zv_1}
  \|y(t)-y^0\|\le c_{y_1}\sqrt{\varepsilon\|y^0-y_0^*\|}+c_{y_2}\varepsilon,
\end{equation}
where
$$c_{y_1} = c_u(M_g+c_g e^{c_h\sqrt{\varepsilon}}M_h\sqrt{\varepsilon})e^{\sqrt{\varepsilon} c_g(1+ \sqrt{\varepsilon}c_h e^{(c_g+c_h)\sqrt{\varepsilon}})},$$
$$c_{y_2} = (M_{g_0}+c_g e^{c_h\sqrt{\varepsilon}}M_{h_0}\sqrt{\varepsilon})e^{\sqrt{\varepsilon} c_g(1+ \sqrt{\varepsilon}c_h e^{(c_g+c_h)\sqrt{\varepsilon}})}.$$

With the obtained estimate, inequality~\eqref{eq8} reads as
\begin{equation}\label{zv_2} \|z(t)-z^0\|\le c_{z_1}\sqrt{\varepsilon\|y^0-y_0^*\|}+\varepsilon c_{z_2},\end{equation}
where $c_{z_1} = e^{c_h\sqrt{\varepsilon}}(M_h c_u+\sqrt{\varepsilon} c_h c_{y_1}),$ $c_{z_2} = e^{c_h\sqrt{\varepsilon}}(M_{h_0}+\sqrt{\varepsilon} c_h c_{y_2}).$

Let us underline that the coefficients $c_{y_1},c_{y_2},c_{z_1}$ and $c_{z_2}$ in estimates~\eqref{zv_1} and~\eqref{zv_2} are monotonically increasing with respect to $\varepsilon$ and $\delta.$

Estimates~\eqref{zv_1} and~\eqref{zv_2} ensure the well-definiteness of the solutions of system~\eqref{eq3} on the time interval $[0,\varepsilon].$ Indeed, estimate ~\eqref{zv_2} means that there is no blow-up of the $z$-component of solutions of system~\eqref{eq3} with initial condition $x^0\in \overline{B_{\delta}(\mathit{Y_0^p})}.$ To show that $y(t)\in B_{\delta}^{'}(y^*(t))$ for all $t\in [0,\varepsilon],$ we exploit the estimate~\eqref{zv_1}:

\begin{equation}\label{eq16}\begin{aligned}
  \|y(t)-y^*(t)\|\le \|y(t)-y^0\|+\|y^*(t)-y_0^*\|+\|y^0-y_0^*\|\\\le c_{y_1}\sqrt{\varepsilon\|y^0-y_0^*\|}+c_{y_2}\varepsilon+L^*\varepsilon+p+ \delta.
\end{aligned}\end{equation}

Thus, to ensure the well-definiteness of the solutions in $D^{'}$ for $t\in[0,\varepsilon],$ it suffices to show that
$$\|y(t)-y^*(t)\|\le {\rm{dist}}(y^*(t),\partial  D^{'})=p+\delta^{'}$$ for each $t\in[0,\varepsilon].$

As $\delta< \delta^{'},$ we may define $\varepsilon_0$ as the positive root of the equation
$$c_{y_1}\sqrt{\varepsilon(p+\delta)}+\varepsilon (c_{y_2}+L^*) = \delta^{'} - \delta,$$
i.e. $$\varepsilon_0 = \left(\sqrt{\left(\frac{c_{y_1}\sqrt{p+\delta}}{2(c_{y_2}+L^*)}\right)^2+\frac{\delta^{'}-\delta}{c_{y_2}+L^*}}-\frac{c_{y_1}\sqrt{p+\delta}}{2(c_{y_2}+L^*)}\right)^2 .$$
Then for any $\varepsilon \in [0,\varepsilon_0],$ the solutions of system~\eqref{eq3} with controls~\eqref{eq5} and initial conditions $x^0\in B_{\delta}(Y_{0}^p)$ are well-defined in $D^{'}$ for all $t\in[0,\varepsilon].$

The next step of the proof is to show that the distance between $y(t)$ and $y^*(t)$ does not increase after the time $t=\varepsilon,$ i.e. $\|y(\varepsilon)-y^*(\varepsilon)\|\le \|y^0-y^*_0\|.$ For this purpose, note that any solution of system~\eqref{eq3} with initial data $x^0 \in B_{\delta}(\mathit {Y_{0}^p} )$ and controls~\eqref{eq5} can be represented by means of the Chen--Fliess type series \cite{La95,ZuSIAM,GZ19a, GZ20_ECC}. For analyzing the value $y(\varepsilon),$ consider the $y$-component of the series expansion, where the term $\varepsilon \mathcal{F}(x^0)a(x^0,y^*_0)$ is added and subtracted:
\begin{equation}\label{eq10}\begin{aligned}
  y(\varepsilon)=y^0&\pm\varepsilon \mathcal{F}(x^0)a(x^0,y^*_0)\\
 & +\int\limits_{0}^{\varepsilon}\left(g_0(t,x)+\sum_{k=1}^{m}g_k(x)u_k^{\varepsilon}(t,x,y^*)\right)dt\\
 =y^0&-\varepsilon\alpha (y^0-y_0^*)+\varepsilon g_0(0,x^0) \\
 &+\sigma_1(\varepsilon,x^0)+r_0(\varepsilon)+r_1(\varepsilon),
\end{aligned}
\end{equation}
where
\begin{equation*}\begin{aligned}
 &\sigma_1(\varepsilon,x^0)=-\varepsilon\mathcal{F}(x^0)a(x^0,y^*_0)+\sum_{k=1}^{m}g_k(x^0)\int\limits_{0}^{\varepsilon}u_k(s_1)ds_1\\
 &\quad+\sum_{k_1,k_2=1}^{m} \mathcal L_{f_{k_2}}g_{k_1}(x^0)\int\limits_{0}^{\varepsilon}\int\limits_{0}^{s_1}u_{k_1}(s_2)u_{k_2}(s_2)ds_2 ds_1,\\
&r_0(\varepsilon)=\int\limits_{0}^{\varepsilon}\int\limits_{0}^{s_1}\Big(\frac{\partial g_0(s_2,x(s_2))}{\partial s_2}+\mathcal L_{f_0}g_0(s_2,x(s_2))\\
&\quad\qquad+\sum_{k=1}^{m}\big(\mathcal L_{f_k}g_0(s_2,x(s_2))u_k(s_2)\\
&\quad\qquad+\mathcal L_{f_0}g_k(s_2,x(s_2))u_k(s_1)\big)\Big)ds_2 ds_1,\\
&r_1(\varepsilon)=\sum_{k_1,k_2=1}^{m}\int\limits_{0}^{t}\int\limits_{0}^{s_1}\int\limits_{0}^{s_2}\Big(\mathcal L_{f_0}\mathcal L_{f_{k_2}}g_{k_1}(x(s_3))\\
&\quad\qquad+\sum_{k_3{=}1}^{m}\mathcal L_{f_{k_3}}\mathcal L_{f_{k_2}}g_{k_1}(x(s_3))u_{k_3}(s_3)\Big)\\
&\quad\qquad\times u_{k_1}(s_1)u_{k_2}(s_2)ds_3ds_2ds_1.
\end{aligned}\end{equation*}
Let us estimate the values of $\|\sigma_1(\varepsilon,x^0)\|$, $\|r_0(\varepsilon)\|$, $\|r_1(\varepsilon)\|.$

Calculating the integrals in $\sigma_1(\varepsilon,x^0)$ according to formula~\eqref{eq5}, we get
\begin{equation*}\begin{aligned}
\sigma_1(\varepsilon,x^0)& = -\varepsilon\mathcal{F}(x^0)a(x^0,y^*_0)+\sum_{k=1}^{m}\varepsilon g_k(x^0) a_k(x^0,y_0^*)\\
  &+\varepsilon \sum_{(k_1,k_2) \in S_2}I_{[n_1\times n]}[f_{k_1},f_{k_2}]a_{k_1 k_2}(x^0,y_0^*)\\
  &+\frac{\varepsilon^2}{2}\sum_{(k_1,k_2) \in S_1} \mathcal L_{f_{k_2}}g_{k_1}(x^0)a_{k_1}(x^0,y_0^*)a_{k_2}(x^0,y_0^*)\\
  &+\frac{\varepsilon^{3/2}}{\sqrt{\pi}}\sum_{k_1  \in S_1}a_{k_1}(x^0,y_0^*)\sum_{k_2 = 1}^m I_{[n_1\times n]}[f_{k_1},f_{k_2}]\\
  &\times\sum_{j:(j,{k_2}) \in S_2} \sqrt{\frac{|a_{j k_2}(x^0,y_0^*)|}{\kappa_{j k_2}}} sign(a_{j k_2}(x^0,y_0^*))
  \\
&  =\frac{\varepsilon^2}{2}\sum_{(k_1,k_2) \in S_1} \mathcal L_{f_{k_2}}g_{k_1}(x^0)a_{k_1}(x^0,y_0^*)a_{k_2}(x^0,y_0^*)\\
  &+\frac{\varepsilon^{3/2}}{\sqrt{\pi}}\sum_{k_1  \in S_1}a_{k_1}(x^0,y_0^*)\sum_{k_2 = 1}^m I_{[n_1\times n]}[f_{k_1},f_{k_2}]\\
  &\times\sum_{j:(j,{k_2}) \in S_2} \sqrt{\frac{|a_{j k_2}(x^0,y_0^*)|}{\kappa_{j k_2}}} sign(a_{j k_2}(x^0,y_0^*)).
\end{aligned}\end{equation*}
Then from A1.3),
\begin{equation*}\begin{aligned}
\|\sigma_1&(\varepsilon,x^0)\|\le\frac{\varepsilon^2 M_{g_2}}{2}\|a(x^0,y^*_0)\|^2\\
&+\frac{2\varepsilon^{\frac{3}{2}}M_{g_2}}{\sqrt{\pi}}\|a(x^0,y^*_0)\|^{\frac{3}{2}}\sum_{j_1=1}^{m}\left(\sum_{(j_2,j_1)\in S_2}\kappa_{j_2j_1}^{-\frac{2}{3}}\right)^{\frac{3}{4}}.\end{aligned}\end{equation*}

By the definition of $a(x^0,y_0^*),$
\begin{equation*}\|\sigma_1(\varepsilon,x^0)\|\le c_{\sigma}\varepsilon^{3/2}\|y_0^*-y^0\|^{3/2},\end{equation*}\
where
\begin{equation*}
\begin{aligned}c_{\sigma} =& \frac{\sqrt{\varepsilon}M_{g_2}\alpha^2 \mu^2}{2}\sqrt{p+\delta}\\
&+\frac{2M_{g_2} (\alpha \mu)^{3/2}}{\sqrt{\pi}}\sum_{j_1=1}^{m}\left(\sum_{(j_2,j_1)\in S_2}\kappa_{j_2j_1}^{-\frac{2}{3}}\right)^{\frac{3}{4}},\end{aligned}\end{equation*}
provided that $x^0\in B_{\delta}(Y_0^p)$

For estimating $\|r_0(\varepsilon)\|$ and $\|r_1(\varepsilon)\|,$ we apply Assumption A1.2):
\begin{equation*}
\|r_0(\varepsilon)\|\le c_{r_0} \varepsilon^\frac{3}{2},\;
\|r_1(\varepsilon)\|\le c_{r_1} \varepsilon^\frac{3}{2} \|y^0-y_0^*\|,
\end{equation*}
where $c_{r_0} = \frac{1}{2} \sqrt{\varepsilon}( \sqrt{\varepsilon}L_{g_{20}}+M_{g_{20}}(\sqrt{\varepsilon}+c_u\sqrt{p+\delta})),\quad\quad$ $c_{r_1} = \frac{1}{6}c_u^2(M_{g_{30}}\sqrt{\varepsilon}+M_{g_3}c_u\sqrt{p+\delta}).$

Let us analyse the value $\|y(\varepsilon)-y^*(\varepsilon)\|.$
From~\eqref{eq10},
\begin{equation*}\begin{aligned}
y(\varepsilon)-y^*(\varepsilon)=&(1-\varepsilon\alpha)(y^0-y^*_0)-(y^*(\varepsilon)-y^*_0)\\
&+\varepsilon g_0(0,x^0) +\sigma_1(\varepsilon,x^0)+r_0(\varepsilon)+r_1(\varepsilon).\end{aligned}\end{equation*}
From A1.2), A1.3), and  the above obtained estimates on $\|\sigma_1(\varepsilon,x^0)\|,$ $\|r_0(\varepsilon)\|,$ $\|r_1(\varepsilon)\|,$
\begin{equation*}
\begin{aligned}
\|y(\varepsilon)&-y^*(\varepsilon)\|\le(1-\varepsilon\alpha)\|y^0-y^*_0\|+L^* \varepsilon +\varepsilon M_{g_0}\\
&+c_{\sigma}\varepsilon^{\frac{3}{2}}\|y^0-y^*_0\|^{\frac{3}{2}}+c_{r_0}\varepsilon^{\frac{3}{2}}+c_{r_1}\varepsilon^{\frac{3}{2}}\|y^0-y^*_0\|\\
&\le \|y^0-y^*_0\|(1-\varepsilon(\alpha -c_{\sigma}\sqrt{\varepsilon(p+\delta)}-c_{r_1}\sqrt{\varepsilon}))\\
&+\varepsilon(L^*+M_{g_0}+c_{r_0}\sqrt{\varepsilon}),\end{aligned}\end{equation*}
provided that $\varepsilon<\varepsilon_1=\frac{1}{\alpha}$ and $x^0\in B_{\delta}(Y_0^p).$

Thus, we achieve the following estimate:
\begin{equation*}
\begin{aligned}
\|y(\varepsilon)-y^*(\varepsilon)\|\le\big(1&-\varepsilon(\alpha-\sqrt{\varepsilon} q)\big)\|y_0-y^*_0\|\\
& +\varepsilon (L^*+M_{g_0}+\sqrt\varepsilon c_{r_0}),
\end{aligned}
\end{equation*}
 where $q=c_{\sigma}\sqrt{p+\delta}+c_{r_1}$  is monotonically increasing with respect to  $\delta.$

  Our next goal is to show the attraction of the $y$-components  of the solution to the $p-$neighborhood of the curve $y^*(t)$.
  Assume that
  $
  \alpha>\frac{\nu(L^*+M_{g_0})}{p}
  $
with some $\nu>1$.

  Using estimate~\eqref{eq16}, we may ensure the following property: if $x^0\in Y_0^{\frac{p}{\nu}}$ then $x(t)\in Y^p_t$ for all $t\in [0,\varepsilon]$ with a small enough $\varepsilon.$ Indeed, let us define $\varepsilon_2$ as the positive root of the equation
  \begin{equation*}
  c_{y_1}\sqrt{\frac{\varepsilon p}{\nu}}+\varepsilon(c_{y_2}+L^*)=\frac{p(\nu-1)}{\nu},
  \end{equation*}
  i.e.
  \begin{equation*}
  \varepsilon_2 = \frac{p\big(\sqrt{c_{y_1}^2+4(\nu-1)(c_{y_2}+L^*)}-c_{y_1}\big)^2}{4\nu(c_{y_2}+L^*)^2}.
  \end{equation*}
  Then estimate~\eqref{eq16} yields
 \begin{equation*}
   \begin{aligned}
\|y(t)-y^*(t)\|\le \|y(t)-y^0\|+\|y^*(t)-y^*_0\|+\|y^0-y^*_0\|\\
  \le c_{y_1}\sqrt{\varepsilon \|y^0-y^*_0\|}+\varepsilon c_{y_2} +\varepsilon L^* +\frac{p}{\nu} \le p,\end{aligned}
 \end{equation*}
  provided that $\|y^0-y_0^*\|\le \frac{p}{\nu}$ and $\varepsilon\le \varepsilon_2.$

Consider two possibilities:
\begin{itemize}
  \item[1.1)] If $x^0 \in Y_0^{\frac{p}{\nu}},$ then, as discussed above, $\|y(t)-y^*(t)\|\le p$ for all $t\in[0,\varepsilon],$ that is $x(t)\in Y_t^p$ for all $t\in[0,\varepsilon]$
  \item[1.2)] If $x^0 \in B_{\delta}(Y_0^p)\backslash Y_0^{\frac{p}{\nu}},$ i.e. $\|y^0-y^*_0\|\ge \frac{p}{\nu},$ then
  \begin{equation*}
      \begin{aligned}
     \|y(\varepsilon)-y^*(\varepsilon)\|&\le\big(1-\varepsilon(\alpha-\sqrt{\varepsilon} q)\big)\|y_0-y^*_0\|\\
& +\frac{\varepsilon (L^*+M_{g_0}+\sqrt\varepsilon c_{r_0})\|y_0-y^*_0\|}{\|y_0-y^*_0\|}\\
\le\bigg(1-\varepsilon\Big(\alpha&-\frac{\nu (L^*+M_{g_0})}{p}\\
&\qquad-\sqrt{\varepsilon}\Big(q+\frac{\nu c_{r_0}}{p}\Big)\Big)\bigg)\|y_0-y^*_0\|\\
\end{aligned}
\end{equation*}
For an arbitrary $\lambda\in\Big(0,\alpha-\frac{\nu(L^*+M_{g_0})}{p}\Big)$,  let us define $\varepsilon_3 = \left(\frac{(\alpha-\lambda) p-\nu(L^*+M_{g_0})}{pq+\nu c_{r_0})}\right)^2$.
Then, for any $\varepsilon\in(0,\varepsilon_3),$
$$
   \|y(\varepsilon)-y^*(\varepsilon)\|\le(1-\varepsilon\lambda)\|y_0-y^*_0\|.
$$
      \end{itemize}
Thus, $\|y(\varepsilon)-y^*(\varepsilon)\|\le \|y^0-y^*_0\|\le p +\delta$ and $x(\varepsilon)\in B_{\delta}(Y_{\varepsilon}^p)\subset D^{'}.$

So we may repeat all the above argumentation for the solutions with the initial condition $x(\varepsilon)$ with the same choice of $\lambda\in\left(0, \alpha-\frac{\nu(L^*+M_{g_0})}{p}\right)$ and $\varepsilon \in \min\big\{\varepsilon_0,\varepsilon_1,\varepsilon_3\big\}.$ This proofs the well-definiteness in $D'$ of the solutions of system~\eqref{eq3} with $x^0 \in B_{\delta}(Y_0^p)$ for $t\in[0,2\varepsilon].$ Besides, we can consider again two cases:
\begin{itemize}
  \item[2.1)] if $\|y(\varepsilon)-y^*(\varepsilon)\|\le\frac{p}{\nu}$ then $\|y(t)-y^*(t)\|\le p$ for $t\in[\varepsilon,2\varepsilon];$
  \item[2.2)] if $\|y(\varepsilon)-y^*(\varepsilon)\|>\frac{p}{\nu}$ then
  $$\|y(2\varepsilon)-y^*(2\varepsilon)\|\le (1-\varepsilon\lambda)\|y(\varepsilon)-y^*(\varepsilon)\|.$$
\end{itemize}
Iterating all above-described steps, we may conclude that the solutions of system~\eqref{eq3} with control~\eqref{eq5} and initial conditions $x^0\in {B_{\delta}(Y_0^p)}$  are well defined in $D'$ for all $t\ge 0.$ Furthermore, if $x^0\in Y_0^{\frac{p}{\nu}}$ then $x(t) \in Y_0^p$ for all $t\ge 0.$ If $x^0\in {B_{\delta}(Y_0^p)}\backslash Y_0^{\frac{p}{\nu}},$ then there exists an $N\in\mathbb{N}$ such that $\|y(t)-y^*(t)\|>p$ for each $t=0, \varepsilon, 2\varepsilon,...,(N-1)\varepsilon$ and $\|y(t)-y^*(t)\|\le \frac{p}{\nu}$ for all $t\in[0,+\infty).$ It remains to describe the behavior of $y(t)$ for an arbitrary $t\in[0,N\varepsilon].$

As follows from the previous argumentation, the following estimate holds for $t=0, \varepsilon, 2\varepsilon,..., (N-1)\varepsilon :$
$$y(j\varepsilon)-y^*(j\varepsilon)\le(1-\varepsilon\lambda)^{j}\|y^0-y^*_0\|\le e^{-\lambda j \varepsilon}\|y^0-y^*_0\|,$$
for $j=0,1,...,N-1.$

For an arbitrary $t \in [0, N\varepsilon],$ denote by $t_{in} = \left[ \frac{t}{\varepsilon}\right]$ the integer part of $\frac{t}{\varepsilon}.$ Notice that $t-t_{in}\varepsilon<\varepsilon,$ then
\begin{equation*}
\begin{aligned}
\|y(t)-y^*(t)\|&\le\|y(t_{in}\varepsilon)-y^*(t_{in}\varepsilon)\|+\|y(t)-y(t_{in}\varepsilon)\|\\
+\|y^*(t)&-y^*(t_{in}\varepsilon)\|\le\sqrt{\|y(t_{in}\varepsilon)-y^*(t_{in}\varepsilon)\|}\\
\times\Big( c_{y_1}\sqrt{\varepsilon}+&\sqrt{\|y(t_{in}\varepsilon)-y^*(t_{in}\varepsilon)\|}\Big)+\varepsilon(c_{y_2}+L^*)\\
&\le \gamma_1(\|y^0-y_0^*\|)e^{-\frac{\lambda t}{2}}+\varepsilon\gamma_2,
\end{aligned}
\end{equation*}
where
$$\gamma_1(\|y^0-y_0^*\|) = e^{\frac{\varepsilon\lambda}{2}}\sqrt{\|y^0-y^*_0\|}\big(c_{y_1}\sqrt{\varepsilon}+e^{\frac{\varepsilon\lambda}{2}}\sqrt{\|y^0-y^*_0\|}\big)$$
is monotonically increasing with respect to ${\|y^0-y^*_0\|,}$ and $\gamma_2 = c_{y_2}+L^*$.

This completes the proof of Theorem~1.

\bibliographystyle{plain}
\bibliography{biblio_nonh}
\end{document}